\RequirePackage{ifpdf}
\ifpdf 
\documentclass[pdftex]{sigma}
\else
\documentclass{sigma}
\fi

\begin{document}

\renewcommand{\PaperNumber}{051}

\FirstPageHeading

\def\fup{\bar f}
\def\xup{\bar x}

\def\aup{\bar\alpha}
\def\ati{\tilde\alpha}
\def\hti{\tilde h}
\def\fti{\tilde f}
\def\fht{\hat f}
\def\hht{\hat h}

\def\prof{\vrule depth 3.5pt width 0pt}
\def\dprof{\vrule depth 5.3pt width 0pt}
\def\sac#1#2#3{\prof{\smash{\hbox{\vtop{\hbox{$#1$}%
\hbox{\raise#2pt\hbox{$\mkern#3mu\hat{}$}}}}}}}
\def\sti#1#2#3{\prof{\smash{\hbox{\vtop{\hbox{$#1$}%
\hbox{\raise#2pt\hbox{$\mkern#3mu\tilde{}$}}}}}}}
\def\sba#1#2#3{\prof{\smash{\hbox{\vtop{\hbox{$#1$}%
\hbox{\raise#2pt\hbox{$\mkern#3mu\bar{}$}}}}}}}

\def\fdo{\sba f 8 4}
\def\xdo{\sba x 8 4}
\def\gdo{\sba g 8 4}
\def\fdoone{\sba {f_1} 8 4}
\def\fdotwo{\sba {f_2} 8 4}
\def\ado{\sba {\alpha} 9 4}
\def\adoone{\sba {\alpha_1} 9 4}

\def\fdt{\sti f 9 4}
\def\hdt{\sti h 9 4}
\def\gdh{\sac g 9 4}

\ShortArticleName{Quantum Painlev\'e Equations: from Continuous to
Discrete}

\ArticleName{Quantum Painlev\'e Equations:\\ from Continuous to
Discrete}

\Author{Hajime NAGOYA~$^\dag$, Basil GRAMMATICOS~$^\ddag$ and
Alfred RAMANI~$^\S$}

\AuthorNameForHeading{H. Nagoya, B. Grammaticos and A. Ramani}

\Address{$^\dag$~Graduate School of Mathematical Sciences, The
University of Tokyo, Japan}
\EmailD{\href{mailto:nagoya@ms.u-tokyo.ac.jp}{nagoya@ms.u-tokyo.ac.jp}}

\Address{$^\ddag$~IMNC, Universit\'e Paris VII \& XI, CNRS, UMR
8165, B\^at. 104, 91406 Orsay, France}
\EmailD{\href{mailto:grammati@paris7.jussieu.fr}{grammati@paris7.jussieu.fr}}

\Address{$^\S$~Centre de Physique Th\'eorique, Ecole
Polytechnique, CNRS, 91128 Palaiseau, France}
\EmailD{\href{mailto:ramani@cpht.polytechnique.fr}{ramani@cpht.polytechnique.fr}}

\ArticleDates{Received March 05, 2008, in f\/inal form May 03,
2008; Published online June 09, 2008}

\Abstract{We examine quantum extensions of the continuous
Painlev\'e equations, expressed as systems of f\/irst-order
dif\/ferential equations for non-commuting objects. We focus on
the Painlev\'e equations II, IV and V. From their auto-B\"acklund
transformations we derive the contiguity relations which we
interpret as the quantum analogues of the discrete Painlev\'e
equations.}


\Classification{34M55; 37K55; 81S99}

 \section{Introduction}

 The word ``quantum'' used in the title needs some qualifying. Historically, this  term was introduced  in relation to the discreteness of the spectrum  of operators like the Hamiltonian and the angular momentum. However, with the blossoming of  quantum mechanics, its use was genera\-li\-sed to the description of non-classical objects, typically non-commuting operators. It is in this sense that we are going to use the term quantum in this  paper: namely, for the description of  equations where the  various components of the dependent variable  do not commute among themselves.

In this paper we shall focus on {\it integrable} equations
involving non-commuting variables. Their interest, in particular,
as far as quantum f\/ield theories are concerned, is obvious. A
large lite\-ra\-ture exists concerning supersymmetric or just
fermionic extensions of integrable evolution equations. It is now
clear that the special properties which characterise integrability
can be extended to the case where the dependent variable involves
fermionic as well as bosonic components.

The quantisation  of low dimensional  integrable Hamiltonian
systems has been the object of  extensive investigations by
Hietarinta and collaborators \cite{Hietarinta,Hietarinta and
Collabolators}.  It was shown, in particular,  that the
preservation of integrability in a  quantum setting often
necessitates  the introduction of purely quantum ({\it i.e.},
explicitly $\hbar$ dependent) terms in both the Hamiltonian and
the invariants. In some cases, these terms  can be absorbed  by
the proper ordering and the introduction of  a~non-f\/lat space
metric. However, it is not clear  whether this suf\/f\/ices  in
all cases.

The ``quantum'' extension  of (continuous) Painlev\'e equations
has been introduced by one of us (HN) in
\cite{Nagoya2004,Nagoya2007}.  For these  paradigmatic integrable
systems, integrability  is not related to the existence of
invariants, but rather to the fact that these nonautonomous
equations can be obtained as the compatibility  condition of a (an
overdetermined) linear system, the Lax pair. (Incidentally,
Novikov~\cite{Novikov} refers to  this compatibility condition as
a quantisation condition for the spectral curve and thus the
deautonomisation process can be considered, in some formal sense,
as a f\/irst kind of quantisation.)

Starting from ``quantum'' continuous Painlev\'e equations we shall
derive  their contiguity relations which, as already  shown in the
commuting case, can be interpreted as discrete Painlev\'e
equations \cite{FGR,JM}. Quantum and  discrete systems  possess a
common character as far as  the phase-space of their dynamics is
concerned. While in the former  case  the surface  of an
elementary cell  is f\/ixed  by the relation $\Delta x\Delta
p=\hbar$ in the latter case  the elementary cell is rigidly
f\/ixed by $\Delta x=\Delta p=1$ (in the appropriate units). The
main  dif\/f\/iculty  in quantising discrete systems,
integrability notwithstanding, lies in the fact that  one must
introduce a commutation rule  consistent with the
evolution~\cite{QN}. This is  a highly  nontrivial problem.  We
have addressed this question  in \cite{GRPN} and \cite{TGRT} where
we have shown that for the  mappings of the QRT~\cite{QRT} family
the following  commutation rule
\begin{gather}
xy=qyx+\lambda x + \mu y +\nu \label{eq:1-1}
\end{gather}
is consistent with the evolution (with the adequate choice of the
parameters).  Since  discrete Painlev\'e equations are
nonautonomous extensions of  the QRT mappings  we expect rule
\eqref{eq:1-1}  to be suf\/f\/icient  for the quantisation of the
cases  we shall consider here.

In what follows, we shall analyse the ``quantum''  forms of
continuous  Painlev\'e equations derived by one of us (HN)  in
\cite{Nagoya2004} and use  their auto-B\"acklund  transformations
so as to  deduce their contiguity relations. Thus  working with
quantum  forms of P$_{\rm II}$, P$_{\rm IV}$ and P$_{\rm V}$  we
will derive quantum forms for  the discrete P$_{\rm I}$, P$_{\rm
II}$, P$_{\rm III}$ and P$_{\rm IV}$.

\section{Non-commuting variables: some basic relations}

Before proceeding to the derivation of quantum discrete Painlev\'e
equations we should present  a~summary of our f\/indings
in~\cite{Nagoya2004}.  Our derivation of quantum continuous
Painlev\'e equations consists in extending the symmetrical form of
Painlev\'e equations
 proposed by Noumi and Yamada \cite{NYHigher} (see also Willox et
 al. \cite{Willox et al}) to non-commuting objects. In the case of the quantum P$\rm _{II}$ equation we introduce  three unknown operators $f_0$, $f_1$, $f_2$ of $t$ and two parameters $\alpha_0$, $\alpha_1$ in the complex number f\/ield $\mathbb{C}$.
The commutation rules are
 \begin{gather}
[f_1,f_0]=2\hbar f_2,\qquad[f_0,f_2]=[f_2,f_1]=\hbar.
\label{eq:2-1}
\end{gather}
(In a more physical parlance we can say that the $\hbar$ appearing
in the commutation relations is the Planck constant).

Generalising  these relations to more objects so as to  describe
higher quantum Painlev\'e equations is straightforward. For a
positive number $l$ greater than $1$, we introduce $l+1$ c-number
parameters $\alpha_i$  $(0\le i\le l)$  and unknown operators
$f_i$ $(0\le i\le l)$ with commutation relations
\begin{gather}\label{eq:2-2}
[f_i,f_{i+1}]=\hbar ,\qquad[f_i,f_j]=0 \quad\rm{otherwise},
\end{gather}
 where the indices $0,1,\ldots,l$ are understood as elements of $\mathbb{Z}/(l+1)\mathbb{Z}$.

In the case of the quantum P$\rm _{II}$ equation,  the evolution
equations with respect to $t$ for the unknown operators $f_i$ are
\begin{gather}\label{eq:2-3}
\partial_t f_0=f_0f_2+f_2f_0+\alpha_0,\qquad  \partial_t f_1=-f_1f_2-f_2f_1+\alpha_1,\qquad  \partial_t f_2=f_1-f_0,
\end{gather}
 and in the case of the quantum P$\rm _{IV}$ equation, the quantum P$\rm _{V}$ equation and higher Painlev\'e equations,
 the evolution equations with respect to $t$ for the unknown operators $f_i$ are: for $l=2n$ ($n\ge 1$),
\begin{gather}\label{eq:2-4}
\partial_t f_i=f_i\left(\sum_{1\le r\le n}f_{i+2r-1}  \right)
-\left( \sum_{1\le r\le n}f_{i+2r}\right)f_i+\alpha_i,
\end{gather}
 and for $l=2n+1$ ($n\ge 1$),
 \begin{gather}
  \partial_t f_i= f_i\left(\sum_{1\le r\le s\le n}f_{i+2r-1}f_{i+2s}\right)
 -\left(\sum_{1\le r\le s\le n}f_{i+2r}f_{i+2s+1}\right)f_i
 \nonumber
\\
\phantom{\partial_t f_i=}{} +\left({k\over 2}-\sum_{1\le r\le
n}\alpha_{i+2r}\right)f_i+\alpha_i\sum_{1\le r\le n}f_{i+2r},
 \label{eq:2-5}
 \end{gather}
 where $k=\alpha_0+\cdots+\alpha_l$. These (non-commutative) dif\/ferential systems are quantization
 of nonlinear ordinary dif\/ferential systems proposed by Noumi and Yamada in \cite{NYHigher},
 and equal to quantization of P$\rm _{IV}$ and P$\rm _{V}$ for $l=2$ and $l=3$, respectively.

 In the introduction we have stressed the necessity for the commutation rule to be consistent with the evolution. All systems here are consistent
with the corresponding evolution $\partial_t$, namely the
evolution $\partial_t$ preserves the commutation relations
\eqref{eq:2-1} for the case of the quantum P$\rm_{II}$ equation
and \eqref{eq:2-2} for the other quantum Painlev\'e equations.

The non-commutative dif\/ferential systems \eqref{eq:2-3},
\eqref{eq:2-4} and \eqref{eq:2-5} admit the af\/f\/ine Weyl group
actions of type $A_l^{(1)}$, respectively, as well as the
classical case. The actions are as follows.
 In the case of the quantum P$\rm _{II}$ equation, we have
  \begin{gather}
s_0(f_ 0)=f_0, \qquad s_0(f_ 1)=f_1-f_2{\alpha_0\over
f_0}-{\alpha_0\over f_0}f_2-{\alpha_0^2\over f_0^2},\qquad
s_0(f_ 2)=f_2+{\alpha_0\over f_0},\nonumber
\\
  s_1(f_ 0)=f_0+f_2{\alpha_1\over f_1}+{\alpha_1\over f_1}f_2-{\alpha_1^2\over f_1^2}, \qquad s_1(f_ 1)=f_1,\qquad   s_1(f_ 2)=f_2-{\alpha_1\over f_1},\nonumber
\\
    s_0(\alpha_ 0)=-\alpha_0, \qquad s_0(\alpha_ 1)=\alpha_1+2\alpha_0,\qquad   s_1 (\alpha_ 0)=\alpha_0+2\alpha_1,\quad s_1(\alpha_ 1)=-\alpha_1,\nonumber
    \\
     \pi(f_ 0)=f_1, \qquad \pi(f_ 1)=f_0,  \qquad \pi(f_ 2)=-f_2, \qquad \pi(\alpha_ 0)=\alpha_1,\qquad   \pi(\alpha_ 1)=\alpha_0.
     \label{eq:2-6}
     \end{gather}
The actions $s_0$, $s_1$ and $\pi$ preserve the commutation
relations \eqref{eq:2-1} and give a representation of the extended
af\/f\/ine Weyl group of type $A_1^{(1)}$, namely, they satisfy
the relations
 \[
s_i^2=1, \qquad \pi^2=1,\qquad \pi s_i=s_{i+1} \pi.
\]
Moreover the actions of $s_0$, $s_1$ and $\pi$ commute with the
dif\/ferentiation $\partial_t$, that is, they are B\"acklund
transformations of the quantum P$\rm _{II}$ equation.

In the case of the quantum P$\rm _{IV}$ equation ($l=2$), the
quantum P$\rm _{V}$ equation ($l=3$) and higher quantum Painlev\'e
equations ($l\ge 4$), for $i,j=0,1,\ldots,l$ ($l\ge 2$) we have
 \begin{gather}\label{eq:2-7}
s_i(f_ j)=f_j+{\alpha_i\over f_i} u_{ij},\qquad s_i(\alpha_
j)=\alpha_j-\alpha_i a_{ij},\qquad \pi(f_ j)=f_{j+1},\quad
\pi(\alpha_ j)=\alpha_{j+1},\
\end{gather}
 where
 \begin{gather}
u_{i,i\pm1}=\pm1,\qquad u_{l0}=1, \qquad u_{0l}=-1, \qquad
u_{ij}=0\quad {\rm otherwise},\nonumber
\\
 a_{ii}=2, \qquad a_{i,i\pm1}=-1,\qquad a_{l0}= a_{0l}=-1, \qquad a_{ij}=0\quad {\rm otherwise}.\nonumber
\end{gather}
The actions of $s_i$ ($i=0,1,\ldots,l$) and $\pi$ preserve the
commutation relations \eqref{eq:2-2} and give representations of
the extended af\/f\/ine Weyl groups of type $A_l^{(1)}$, namely,
they satisfy the relations
\[
s_i^2=1,\qquad (s_is_{i+1})^3=1,\qquad s_is_j=s_js_i\quad (j\neq
i\pm 1),\qquad \pi^{l+1}=1,\quad \pi s_i=s_{i+1}\pi.
\]
Moreover the actions of $s_i$ ($i=0,1,\ldots,l$) and $\pi$ commute
with the dif\/ferentiation $\partial_t$, that is, they are
B\"acklund transformations of the corresponding non-commutative
dif\/ferential systems.

We stress that the auto-B\"acklund transformations preserve the
commutation relations in each case so that a discrete evolution
which is constructed from these auto-B\"acklund transformations
also preserves the commutation relations.

We remark that each dif\/ferential system \eqref{eq:2-3},
\eqref{eq:2-4} or \eqref{eq:2-5} has  relations
\[
\partial_t(f_0+f_1+f_2^2)=k
\]
for the quantum P$\rm _{II}$ case,
\[
\partial_t\left(\sum_{r=0}^lf_r\right)=k
\]
for the $l=2n$ ($n\ge 1$) case and
\[
\partial_t\left(\sum_{r=0}^lf_{2r}\right)={k\over 2}\sum_{r=0}^lf_{2r},\qquad
\partial_t\left(\sum_{r=0}^lf_{2r+1}\right)={k\over 2}\sum_{r=0}^lf_{2r+1}
\]
for the $l=2n+1$ ($n\ge 1$) case, where
$k=\alpha_0+\alpha_1+\cdots+\alpha_l$. For simplicity, we
normalize $k=1$ in the following.

\section{The continuous quantum Painlev\'e II\\ and the related discrete equation}

As explained in Section~2 the Painlev\'e II case can be obtained
with three non-commuting  objects $f_0$, $f_1$, $f_2$ as described
in \eqref{eq:2-1} and parameters $\alpha_0$, $\alpha_1$. From the
dynamical equations
\begin{gather}
f_0'=f_0f_2+f_2f_0+\alpha_0,\nonumber
\\
f_1'=-f_1f_2-f_2f_1+\alpha_1,\nonumber
\\
f_2'=f_1-f_0,\nonumber
\end{gather}
(where the prime $'$ denotes dif\/ferentiation with respect to
$t$), eliminating $f_0$, $f_1$ we obtain the equation for the
quantum P$_{\rm II}$
\begin{gather*}
f_2''=2f_2^3-tf_2+\alpha_1-\alpha_0,
\end{gather*}
which has the same expression as the commutative P$_{\rm II}$. On
the other hand, if we eliminate~$f_0$,~$f_2$ we f\/ind the quantum
version of P$_{34}$. We obtain
\begin{gather*}
f_1''={1\over2}f_1'f_1^{-1}f_1'-4f_1^2+2tf_1-{1\over2}(\alpha_1^2-\hbar^2)f_1^{-1}.
\end{gather*}
The dif\/ference of this ``quantum'' version with the commutative
P$_{34}$ is clearly seen in the term quadratic in the f\/irst
derivative (which would have been $f_1'^2/f_1$ in the commutative
case) but also in the last term. As a matter of fact, the
coef\/f\/icient of the term proportional to $1/f_1$ is a~perfect
square in the commutative case. The appearance of the $-\hbar^2$
in the $\alpha_1^2-\hbar^2$ coef\/f\/icient is a consequence of
the non-commutative character of the $f_i$'s.

In order to derive the discrete equation obtained as a contiguity
relation of the solutions of the quantum P$_{\rm II}$ we shall use
the relations \eqref{eq:2-6} presented in the previous section. We
def\/ine an evolution in the parameter space of P$_{\rm II}$ by
$\fup\equiv s_1\pi f$ and as a consequence the reverse evolution
is $\fdo\equiv\pi s_1f$. Using the relations \eqref{eq:2-6} we
f\/ind
\begin{gather}\label{eq:3-1}
\fup_2+f_2=\alpha_1f_1^{-1}
\end{gather}
and
\begin{gather}\label{eq:3-2}
\fdoone+f_1=t-f_2^2.
\end{gather}
In the same way we compute the ef\/fect of the transformations on
the parameters $\alpha$. We f\/ind $\aup_1=\alpha_1+1$ (and
similarly $\adoone=\alpha_1-1$). Thus applying $n$ times the
transformation $s_1\pi$ on $\alpha_1$, we f\/ind that $\alpha_1$
becomes $\alpha_1+n$. We eliminate $f_1$ between the two equations
\eqref{eq:3-1} and \eqref{eq:3-2} and obtain an equation governing
the evolution in the parameter space. We f\/ind
\begin{gather*}
\alpha_1(\fup_2+f_2)^{-1}+\adoone(\fdotwo+f_2)^{-1}=t-f_2^2.
\end{gather*}
This is the quantum version of the equation known, in the
commutative case, as the alternate discrete Painlev\'e~I.

\section{The continuous quantum Painlev\'e IV\\ and the related discrete equation}

We turn now to the case of the quantum P$_{\rm IV}$ which is
obtained from the equations presented in Section~2 for $l=2$.
Again we have three dependent variables $f_0$, $f_1$, $f_2$. The
quantum continuous P$_{\rm IV}$ equation is
\begin{gather}
f_0'=f_0f_1-f_2f_0+\alpha_0,\nonumber
\\
f_1'=f_1f_2-f_0f_1+\alpha_1,\nonumber
\\
f_2'=f_2f_0-f_1f_2+\alpha_2, \nonumber
\end{gather}
which is very similar to the one obtained in the commutative case.
It is interesting to eliminate $f_0$ and $f_2$ and give the
equation for $f_1$. After a somewhat lengthy calculation we f\/ind
\begin{gather*}
f_1''={1\over2}f_1'f_1^{-1}f_1'+{3\over2}f_1^3-2tf_1^2+\left({t^2\over2}+\alpha_2-\alpha_0\right)f_1-{1\over2}(\alpha_1^2-\hbar^2)f_1^{-1},
\end{gather*}
where we have used the identity $\sum_{i}f_i=t$. Again we remark
that this equation, with respect to the commutative P$_{\rm IV}$,
contains an explicit quantum correction.

Before proceeding to the construction of the discrete system
related to this equation we derive some auxiliary results.
Starting from the action of $s_i$ on $f_j$, which from
\eqref{eq:2-7} is just $s_i(f_ j)=f_j+\alpha_i f_i^{-1} u_{ij}$ we
f\/ind $s_i(\sum_{j}f_ j)=\sum_{j}f_ j$ because $\sum_{j}
u_{ij}=0$. Similarly $\pi(\sum_{j}f_ j)=\sum_{j}f_ j$. Thus
$\sum_{j}f_ j$ is conserved under any combination of the
transformations $\pi$ and $s_i$ of \eqref{eq:2-7}.

We can now def\/ine the evolution in the parameter space just as
we did in the case of the quantum P$_{\rm II}$. We put $\fup\equiv
s_1s_0\pi f$ and for  the reverse evolution is $\fdo\equiv\pi^{-1}
s_0s_1f$. Using the relations  \eqref{eq:2-7} we f\/ind
\begin{gather}\label{eq:4-1}
\fup_2+f_1+f_2=t-{\alpha_1\over f_1}
\end{gather}
and similarly
\begin{gather}\label{eq:4-2}
f_1+f_2+\fdoone=t+{\alpha_2\over f_2}.
\end{gather}
We now study the ef\/fect of the up-shift operator on $\alpha_1$
and $-\alpha_2$ which play the role of the independent variable.
(The minus sign in front of $\alpha_2$ guarantees that the two
equations have the same form). Using the relations in
\eqref{eq:2-7} and the fact that the sum of the $\alpha$'s is
constant we f\/ind $\aup_1=\alpha_1+1$ and
$(-\aup_2)=(-\alpha_2)+1$. Thus incrementing the independent
variable under repeated applications of the up-shift operator
leads to a linear dependence on the number of iterations. Still,
since the starting point is dif\/ferent we have two free
parameters. Thus the system \eqref{eq:4-1}, \eqref{eq:4-2} is
exactly the quantum analogue of the equation known (in the
commutative case) as the ``asymmetric discrete Painlev\'e I'',
which is, in fact, a discrete form of P$_{\rm II}$ \cite{ADPI}.

\section{The continuous quantum Painlev\'e V\\ and the related discrete systems}

Finally we examine the case of the quantum P$_{\rm V}$,
corresponding to the case $l=3$ in Section~2. The quantum
continuous P$_{\rm V}$ equation is
\begin{gather}
f_0'=f_0f_1f_2-f_2f_3f_0+\left({1\over
2}-\alpha_2\right)f_0+\alpha_0f_2,\nonumber
\\
f_1'=f_1f_2f_3-f_3f_0f_1+\left({1\over
2}-\alpha_3\right)f_1+\alpha_1f_3,\nonumber
\\
f_2'=f_2f_3f_0-f_0f_1f_2+\left({1\over
2}-\alpha_0\right)f_2+\alpha_2f_0,\nonumber
\\
f_3'=f_3f_0f_1-f_1f_2f_3+\left({1\over
2}-\alpha_1\right)f_3+\alpha_3f_1.\nonumber
\end{gather}
We remark that $f_0'+f_2'=(f_0+f_2)/2$ and similarly
$f_1'+f_3'=(f_1+f_3)/2$. Thus two of the variables can be easily
eliminated by introducing explicitly the time variable through the
exponential $e^{t/2}$. Just as in the case of P$_{\rm III}$ and
P$_{\rm IV}$ we can eliminate one further variable and obtain the
quantum form of P$_{\rm V}$ expressed in terms of a single
variable. In the present case it is more convenient to introduce
an auxiliary variable $w=1-e^{t/2}/f_0$. We obtain thus for $w$
the equation
\begin{gather}
w''=w'\left({1\over
w-1}+{1\over2w}\right)w'+(w-1)^2\left({\alpha_0^2-\hbar^2\over2}w+{\hbar^2-\alpha_2^2\over2w}\right)\nonumber
\\
\phantom{w''=}{}
+(\alpha_3-\alpha_1)e^tw+e^{2t}{w(w+1)\over2(1-w)}.\label{eq:5-1}
\end{gather}
As in the previous cases, the equation \eqref{eq:5-1} contains an
explicit quantum correction, with respect to the commutative
P$_{\rm V}$, as well as a symmetrisation of the term quadratic in
the f\/irst derivative.

Just as in the case of commutative P$_{\rm V}$ we can def\/ine
several dif\/ferent evolutions in the parameter space giving rise
to dif\/ferent discrete equations. For the f\/irst equation we
introduce the f\/irst up-shift operator $R=\pi s_3s_2s_1$ from
\eqref{eq:2-7} for $l=3$. If we def\/ine $x=f_0+f_2$ and
$y=f_1+f_3$ a careful application of the rules \eqref{eq:2-7}
shows that $Rx=y$ and $Ry=x$. Next we introduce an auxiliary
variable $g$ def\/ined by
$g=f_3-\alpha_0f_0^{-1}=y-f_1-\alpha_0f_0^{-1}$ and we seek an
equation in terms of the variables $f_1$, $f_0$,  $g$. Calling
$\fup\equiv Rf$, we f\/ind
\begin{gather}\label{eq:5-2}
\fup_1+f_0=x-{\alpha_0+\alpha_3\over g}
\end{gather}
complemented by the equation coming from the def\/inition of $g$
\begin{gather}\label{eq:5-3}
f_1+g=y-{\alpha_0\over f_0}
\end{gather}
and f\/inally
\begin{gather}\label{eq:5-4}
\gdo+f_0=x+{\alpha_1\over f_1}.
\end{gather}
In order to introduce the proper independent variable we consider
the action of $R$ on the $\alpha$'s: $\aup_0=\alpha_0+1$,
$\aup_1=\alpha_1-1$, $\aup_2=\alpha_2$ and $\aup_3=\alpha_3$. Thus
we can choose $z=\alpha_0$ (which grows linearly with the
successive applications of $R$). The system \eqref{eq:5-2},
\eqref{eq:5-3} and \eqref{eq:5-4} is the quantum analogue of the
asymmetric, ternary, discrete Painlev\'e I, which, as we have
shown in~\cite{dressing} is a~discrete form of the Painlev\'e IV
equation. (At this point we should point out that the alternating
constants~$x$,~$y$ do not introduce  an extra degree of freedom.
As a matter of fact by choosing an appropriate gauge of the
dependent variables and rescaling of the independent ones we can
bring these constants to any non-zero value).

For the second equation we introduce a new up-shift operator
$T=s_1\pi s_3s_2$, the action of which is again obtained with the
help of \eqref{eq:2-7}. The variables $x$ and $y$ are def\/ined in
the same way as in the previous paragraph and again we have $Tx=y$
and $Ty=x$. An auxiliary variable is necessary in this case also
and thus we introduce $h=f_3+\alpha_2f_2^{-1}$. We denote the
action of $T$ by a tilde: $\fti\equiv Tf$. We seek an equation for
$f_1$, $f_2$, $h$. We f\/ind
\begin{gather}\label{eq:5-5}
\fdt_1+f_2=x+{\alpha_2+\alpha_3\over h}.
\end{gather}
The def\/inition of $h$ implies
\begin{gather}\label{eq:5-6}
h+f_1=y+{\alpha_2\over f_2}
\end{gather}
and f\/inally we have
\begin{gather}\label{eq:5-7}
f_2+\hti=x-{\alpha_1\over f_1}.
\end{gather}
The action of $T$ on the $\alpha$'s is $\ati_0=\alpha_0$,
$\ati_1=\alpha_1+1$, $\ati_2=\alpha_2-1$ and $\ati_3=\alpha_3$.

Combining the two systems above we can obtain a nicer, and more
familiar, form. First,  comparing \eqref{eq:5-4} and
\eqref{eq:5-7}, we f\/ind that
$\gdo=f_2+\alpha_1f_1^{-1}=x-\hti=T(y-h)$, or equivalently
$g=RT(y-h)$. We are thus led to introduce the operator $S=RT$, and
we denote its action by a ``hat'' accent $\fht=Sf$. Subtracting
\eqref{eq:5-3} from \eqref{eq:5-6} and eliminating $g$ through the
relation $g=y-\hht$, we obtain
\begin{gather}\label{eq:5-8}
h+\hht=y+{\alpha_0\over f_0}-{\alpha_2\over f_0-x}.
\end{gather}
For the second equation we start from \eqref{eq:5-5} and apply the
operator $RT$ to it. We f\/ind
\begin{gather}\label{eq:5-9}
\fup_1+\fht_2=x-{\alpha_0+\alpha_1\over \hht}
\end{gather}
where we have used the fact that
$RT(\alpha_2+\alpha_3)=-(\alpha_0+\alpha_1)$. Next we subtract
\eqref{eq:5-9} from \eqref{eq:5-2}, use the def\/inition of $x$ in
order to eliminate $f_2$ and f\/ind f\/inally
\begin{gather}\label{eq:5-10}
f_0+\fht_0=x+{\alpha_0+\alpha_3\over
\hht-y}+{\alpha_0+\alpha_1\over \hht}.
\end{gather}
A careful application of the rules \eqref{eq:2-7} shows that the
action of the operator $RT$ on each of the numerators of the
fractions in the right hand side of \eqref{eq:5-8} and
\eqref{eq:5-10} results in  an increase by exactly 1 and thus an
independent variable linear in the number of iterations of $RT$
can be introduced. As in the case of the f\/irst system
\eqref{eq:5-2}, \eqref{eq:5-3} and \eqref{eq:5-4} presented in
this section the alternating constants $x$, $y$ do not introduce
an extra degree of freedom. A better choice of the dependent
variables would be $h\to h-y/2$ and $f_0\to f_0-x/2$. Moreover, by
choosing an appropriate gauge of the dependent variables and
rescaling of the independent ones we can bring these constants to
any non-zero value, for instance $x=y=2$. The system now becomes
\begin{gather}
h+\hht={\alpha_0\over f_0+1}-{\alpha_2\over f_0-1},
\label{eq:5-11}
\\
f_0+\fht_0={\alpha_0+\alpha_3\over \hht-1}+{\alpha_0+\alpha_1\over
\hht+1}.\label{eq:5-12}
\end{gather}
Under this form one recognizes immediately the structure of the
``asymmetric discrete Painle\-v\'e~II'' equation, introduced in
\cite{PIII} and which is a discrete analogue of P$_{\rm III}$.
Thus \eqref{eq:5-11} and \eqref{eq:5-12} constitute the quantum
extension of the latter.

 \section{Conclusions}

 In this paper, we have analysed the ``quantum''  forms of  Painlev\'e equations derived by one of us (HN)  in \cite{Nagoya2004}. The derivation consists in extending the symmetrical form of Painlev\'e equations proposed by Noumi and Yamada to non-commuting variables. We have focused here on P$_{\rm II}$, P$_{\rm IV}$ and P$_{\rm V}$ and derived their more ``familiar'' forms expressed in terms of a single variable. The non-commutative character manifests itself in the fact that the dependent function does not commute with its f\/irst derivative. As a consequence (and despite the fact that a symmetrised form of the term quadratic in the f\/irst derivative is used) explicit quantum corrective terms appear in the equation, proportional to the square of the Planck constant. Using  the auto-B\"acklund  transformations  of the continuous Painlev\'e equations we  derive their contiguity relations which are just the quantum forms for  the discrete P$_{\rm I}$, P$_{\rm II}$, P$_{\rm III}$ and P$_{\rm IV}$.

\pdfbookmark[1]{References}{ref}
\LastPageEnding
\end{document}